\documentclass[CJK,11pt]{article}
\usepackage{amscd}
\usepackage{mathrsfs}
\usepackage{amsfonts}
\usepackage{amsmath,amscd,stmaryrd,amssymb,array, amsfonts,amsmath,amssymb,graphicx}

\numberwithin{figure}{section}
 \numberwithin{equation}{section}

\newtheorem{theorem}{Theorem}[section]
\newtheorem{proposition}[theorem]{Proposition}
\newtheorem{definition}[theorem]{Definition}

\newtheorem{lemma}[theorem]{Lemma}
\newtheorem{remark}[theorem]{Remark}

\newcommand{\cB}{\mb{B}}

\newcommand{\cM}{{\mathcal M}}

\newcommand{\cR}{{\mathcal R}}

\def\be{\begin{equation}}
\def\ee{\end{equation}}
\def\ba{\begin{array}}
\def\ea{\end{array}}
\def\benu{\begin{enumerate}}
\def\eenu{\end{enumerate}}
\def\bt{\begin{theorem}}
\def\et{\end{theorem}}
\def\bl{\begin{lemma}}
\def\el{\end{lemma}}
\def\br{\begin{remark}}
\def\er{\end{remark}}
\def\ss{\subset\subset}

\def\b{\beta}
 \def\de{\delta} \def\pa{\partial} 
\def\lam{\lambda} 
\def\ve{\varepsilon} \def\sig{\sigma}
\def\w{\omega}\def\W{\Omega}
\def\a{{\mathscr A}}

\def\.{\cdot}

\def\R{\mathbb{R}}

\def\A{\forall}
\def\ol{\overline}

\def\Cap{\bigcap}\def\Cup{\bigcup}
\def\ln{\mbox{ln\,}}
\def\ra{\rightarrow}
\def\stac{\stackrel}
\def\~{\stac{\sim}}
\def\8{\infty}
\def\X{\times}
\def\({\left(}
\def\){\right)}

\def\E{\exists}
\def\mb{\mbox}
\def\ep{\emptyset}
\def\-{\setminus}

\def\Hs{\hspace{1cm}}\def\hs{\hspace{0.5cm}}
\def\Vs{\vskip8pt}\def\vs{\vskip4pt}

\def\({\left(}\def\){\right)}

\begin{document}

\begin{center}
{\bf\Large  Morse Theory of Attractors via Lyapunov Functions}
\end{center}
\vskip40pt
\centerline{Desheng  Li }  %\vskip20pt
\begin{center}
{\footnotesize
{Department of Mathematics, School of Science, Tianjin University\\
     Tianjin 300072,  China\\
{\em E-mail:}  lidsmath@tju.edu.cn, \,\,lidsmath@hotmail.com}}
\end{center}

\vskip50pt

\begin{minipage}{13.5cm}
\centerline{\large\bf Abstract} \vskip20pt This paper is concerned
with the Morse  theory of attractors for semiflows on complete
metric spaces. First, we construct global Morse-Lyapunov functions
for Morse decompositions of attractors. Then we extend some well
known deformation results in the critical-point theory to
Morse-Lyapunov functions which are only continuous. Based on these
works, we finally introduce the concept of
 critical groups for Morse sets and
establish Morse inequalities and Morse equations for attractors.


\Vs\Vs {\bf Running Head:}  Morse Theory of Attractors.

\Vs\Vs {\bf Date:}   April 20, 2009.

\end{minipage}

\newpage

\section{Introduction}

This paper is concerned with the Morse theory of attractors for
semiflows on complete (not necessarily locally compact) metric
spaces.

The attractors of a semiflow are of crucial importance, this is
because that much of the longtime dynamics is represented by the
dynamics on and near the attractors. Of special interest is the
global attractor. Of course not every semiflow has a global
attractor. However many dissipative systems do, and when this global
attractor exists, it is the depository of ``\,all\,'' the longtime
dynamics of the given system.

The existence of attractors (especially  attractors for infinite
dimensional systems) has been extensively studied in the past
decades for both autonomous and nonautonomous systems; see, e.g.,
\cite{Bab,CV,Cheban,Hale,Har,Lady,Ma,Rob,Sell,Sch,Tem} and
\cite{Vishik}.
 In most examples coming from different evolution problems, one can
 give an estimate on the
  Hausdorff (or fractional, or informational)
dimension of an attractor. It can even be proved for an
infinite-dimensional system that the global attractor is actually
contained in a finite-dimensional manifold. Indeed, should such a
manifold exist, the analysis of the longtime  behavior of the
original system will be reduced to the finite-dimensional case. This
question is addressed by the theory of inertial manifolds; see
\cite{Const,Tem} for details. On the other hand, even in the
finite-dimensional situation,  if the attractor of a dissipative
system has pathological geometry, then the asymptotic behavior of
the system can  be very complicated, as in the case of the Lorenz
one. To have a better understanding on the dynamics of a semiflow it
is thus necessary for us to take a look at the structure of the flow
inside its attractors.

Morse decomposition describes geometric and topological structures
of flows  inside attractors and invariant sets. It has important
applications not only in  dynamical systems theory itself, but also
in many other different areas; see, e.g.,
\cite{chen,CK,Hale,Kappos,Lisiam,Lidcds,Mis,PK,Ryba}. It is well
known that an attractor of a gradient system with finite number of
equilibria has a natural Morse decomposition with each Morse set
being precisely an equilibrium.  The general theory concerning
 Morse decompositions for autonomous systems can be found in
the original work of Conley \cite{Conley} (see also \cite{Akin,
Mis,Ryba}). Extensions to non-autonomous systems can be found in
\cite{Ras}, and to random dynamical systems in \cite{Cra, Liu} etc.

In this present work we want to go a further step towards Morse
decomposition theory of attractors for semidynamical systems on
complete metric spaces. Our main purpose  is to develop an easy
approach which allows us to work with  Morse theory of attractors
completely in the framework  of the classical Morse theory for
smooth functionals. First, we construct some nice global
Morse-Lyapunov functions for  Morse decompositions of  attractors,
where {\em global} means that the functions are defined on the whole
attraction basins of attractors in comparison with those in the
literature which only have definitions on the attractors themselves.
More precisely, let $X$ be a complete metric space, and let $\a$ be
an attractor of the semiflow $S(t)$  on $X$ with attraction basin
$\W=\W(\a)$ and Morse decomposition $\cM=\{M_1,\cdots,M_l\}.$ We
will construct  a Lyapunov function $V\in C(\W)$ such that
\begin{enumerate}
 \item[(1)] $\lim_{x\ra\pa\W}V(x)=+\8$;\,\, (radial unboundedness)
\item[(2)] $V$ is constant on each Morse set with
\be\label{e:0.4} V(M_1)<V(M_2)<\cdots<V(M_l); \ee
\item[(3)] there is a nonnegative function $v\in C(\W)$ with
$v(x)>0$ (in $\W\setminus \cR$) such that $$D^+V(x)<-v(x),\Hs x\in
\W\setminus\cR,$$ where $D^+V(x)$ is the Dini derivative of $V$
along the flow, and $$\cR=\Cup_{1\leq k\leq l}M_k.$$
\end{enumerate} In addition, if $X=\R^m$, then for any given $k\geq 1$  we can construct $V$ so that it belongs to  $C^k(\W)$.

 Since Morse decomposition theory plays an important role in the dynamical theory of nonlinear control systems \cite{chen,CK,Kappos},
 such global Morse-Lyapunov functions can be expected to have significant applications  in the design of feedback controls for nonlinear systems,
 therefore they are also of independent interest.

Second, we are interested in the deformation of level sets of
Morse-Lyapunov functions. As one of our main concerns here, we will
try to extend some classical deformation results in the
critical-point theory on smooth functionals to Morse-Lyapunov
functions of semiflows. Specifically, let $V$ be a strict
Morse-Lyapunov function of $\cM$ on $\W$. We will show that if
$V^{-1}([a,b])$ contains no Morse sets, then $V_a$ is a strong
deformation retract of $V_b$, where $V_R$ denotes the level set of
$V$,
$$V_R=\{x\in \W|\,\,V(x)\leq R\}.
$$
If $M_k$ consists of exactly one equilibrium of the semiflow, we
further prove that $V_c$ is a strong deformation retract of $V_b$,
where $c=V(M_k)$, and $b>c$ is any number such that $V^{-1}((c,b])$
contains no Morse sets.

Deformation results are of crucial importance in the critical-point
theory and variational methods \cite{chang2,Stru}. Note that our
results here only require continuity of Morse-Lyapunov functions,
therefore we
 hope that some ideas and techniques used here will be helpful in
developing the  critical-point theory for nonsmooth variational
functionals.

Third, we try to develop Morse theory of attractors completely  in
the framework of the classical Morse theory by using global
Morse-Lyapunov functions and deformation lemmas. First, we introduce
the concept of  critical groups  for Morse sets. Let $V$ be a strict
Morse-Lyapunov function (i.e., $V$ satisfies (\ref{e:0.4})\,), and
let $a<c<b$ be such that $c=V(M_k)$ is the unique generalized
critical value of $V$ in $[a,b]$. (Generalized critical values are
the values of $V$ on Morse sets.) Then we define the { critical
group} $C_*(M_k)$ of $M_k$  to be the sequence of singular homology
groups:
$$
C_q(M_k)=H_q\(V_b,\,V_a\),\Hs q=0,1,\cdots,
$$
where $H_*$ denotes the usual relative  singular homology theory of
space pairs. We show that $C_*(M_k)$ is independent of the choice of
numbers $a$ and $b$ as well as the Morse-Lyapunov functions. As a
matter of fact, we actually  prove that \be\label{e:0.1}
C_q(M_k)\cong H_q\(W_k,\,W_{k-1}\) \ee for any positively invariant
neighborhoods $W_k\subset \W(A_k)$ of $A_k$ and $W_{k-1}\subset
\W(A_{k-1})$ of $A_{k-1}$, where
$$\ep=A_0\nsubseteq A_1\nsubseteq\cdots\nsubseteq A_l=\a$$ is the Morse filtration of $\cM$, and $\W(A_k)$ denotes the
attraction basin of  $A_k$. (Lemma \ref{l:2.9} below implies
$\W(A_k)$ is an open neighborhood of $A_k$.) This makes the
computation of the critical groups easier and more flexible.

Then we establish Morse inequalities and equations for attractors.
Let$$ m_q=\sum_{1\leq k\leq l}\mb{rank}\,C_q\(M_k\)$$ be the $q$-th
Morse type number of $\cM$. We prove that if all the critical groups
are of finite rank, then for any  $q\geq 0$ we have\be\label{e:0.2}
m_q-m_{q-1}+\cdots+(-1)^qm_0\geq
\beta_q-\beta_{q-1}+\cdots+(-1)^q\beta_0,\ee where $\b_q$ is the
$q$-th Betti number of the attraction basin $\W$. In addition,
\be\label{e:0.3}
\sum_{q=0}^\8(-1)^qm_q=\sum_{q=0}^\8(-1)^q\beta_q=\chi(\W),\ee
provided that the left-hand side of the above equation is
convergent, where $\chi(\W)$ is the Euler number of $\W$.

If  $\a$ is the global
 attractor, then $\b_q$
 is precisely the $q$-th Betti number of the phase space  $X$. In
 such a case it is interesting to note that the right-hand sides of (\ref{e:0.2})
and (\ref{e:0.3})  are independent of the concrete systems and
attractors. This suggests  that the Euler number $\chi(X)$ is
somewhat an invariant for dissipative systems on $X$.

A very particular but important case is that $X=M^n$ is an
$n$-dimensional compact $C^1$-manifold, in which all the critical
groups are of finite rank. Let $\cM=\{M_1,\cdots,M_l\}$ be a Morse
decomposition of $M^n$ induced by any semiflow $S(t)$ on $M^n$. Then
(\ref{e:0.2}) and (\ref{e:0.3}) reduce to
$$
m_0\geq \beta_0,$$
$$
m_1-m_0\geq \beta_1-\beta_0,
$$
$$
\cdots\cdots
$$
$$
m_n-m_{n-1}+\cdots+(-1)^nm_0=\beta_n-\beta_{n-1}+\cdots+(-1)^n\beta_0=\chi(M^n).
$$

The quotient critical group $\~C_*(M_k)$ of Morse sets and Morse
inequalities and equations  are also addressed in the quotient phase
space $\~X$.

Morse theory of attractors can be established by using Conley index
theory; see \cite{Mis,Ryba,RZ}. Because of  the potentially
pathological geometry of attractors and Morse sets, in general the
computation of homotopy types of these sets may be  somewhat a
difficult problem. Recently Kapitanski and Rodnianski developed an
alternative approach to construct Morse theory of attractors by
making use of the shape theory and $\check{\mb{C}}$ech homologies
\cite{Kap}. In contrast, the approach here seems to be more direct
and simple, and  is easier to be handled.

\vs This paper is organized as follows.  Section 2 is devoted to
some preliminary works, and Section 3 is concerned with Lyapunov
functions of attractors. In Section 4 we construct strict
Morse-Lyapunov functions for Morse decompositions, and in Section 5
we prove two deformation lemmas. Section 6 consists of the
discussions on critical groups of Morse sets and Morse inequalities
as well as Morse equations of attractors.

\section{Preliminaries }

In this section we recall some basic definitions and facts in the
theory of dynamical systems for semiflows on complete metric spaces.

\subsection{Semiflows}
Let $X$ be  a complete  metric space with metric $d(\cdot,\cdot)$.
For convenience we will always identify a singleton $\{x\}$ with the
 point $x\in X$.

 Let $A$ and $B$ be nonempty subsets of $X$. Define the { distance} $\mb{d}(A,B)$ between $A$ and $B$ as
$$\mb{d}(A,B)= \inf_{x\in A,\,y\in B} d(x,y).
$$
 The
closure and interior of $A$ are denoted by $\ol A$ and int$\,A$. A
subset $U$ of $X$ is called a { neighborhood} of $A$, this means
that $\ol A\subset \mbox{int}\,U$. The {boundary} and
{$\ve$-neighborhood} of $A$ are defined, respectively, as $$\pa
A=\ol A\setminus \mb{int}\,A,\Hs \cB(A,\ve)=\{y\in
X:\,d(y,A)<\ve\}.$$

\begin{definition}
A {semiflow} $($semidynamical system$)$ on $X$ is a continuous
mapping $S:\R^+\X X\ra X$ that satisfies $$ S(0,x)=x,\Hs
S(t+s,x)=S\(t,\,S(s,x)\)$$ for all $x\in X$ and $t,s\geq0.$
\end{definition}

We usually write $S(t,x)$ as $S(t)x$. Therefore a semiflow $S$ can
 be viewed as a family of operators $\{S(t)\}_{t\geq0}$
satisfying:
$$
S(0)=\mb{id}_X,\Hs S(t+s)=S(t)S(s)\,\,\,(\A\,t,s\geq0).$$

\br\label{r:2.0} Since $S(t)x$ is continuous in $(t,x)$, it is
uniformly continuous in $(t,x)$ on any compact set $[0,T]\X K\subset
\R^+\X X$. Consequently we see that if $x_n\ra x_0$, then for any
$T>0$, $S(t)x_n$ converges to $S(t)x_0$ uniformly with respect to
$t$ on any compact interval $[0,T]$. \er

From now on  we will always assume that there has been given a
semidynamical system $S(t)$ on $X$; moreover, we assume $S(t)$ is
{\bf asymptotically  compact}, that is, $S(t)$ satisfies the
following assumption:
\begin{enumerate}
\item[({AC})] For any bounded  sequence $x_n\in X$ and $t_n\ra+\8$, if the
sequence $S(t_n)x_n$ is bounded, then it has a convergent
subsequence.
\end{enumerate}

The  asymptotic compactness property (AC) is fulfilled by a large
number of infinite dimensional semiflows generated by PDEs in
applications \cite{Sell}.

 Let $I\subset\R$ be an interval. A { trajectory} $\gamma$ of $S(t)$ on
$I$ is a mapping $\gamma:I\rightarrow X$ satisfying $\gamma(t)=
S(t-s)\gamma(s)$ for any $s,t\in I$ with $s\leq t$. In case
$I=(-\8,\,+\8)$, we will simply say that $\gamma$ is a { complete
trajectory}. A complete trajectory $\gamma$ through $x\in X$ means a
complete trajectory with $\gamma(0)=x$.

\vs

The following basic fact can be deduced by using the asymptotic
compactness property of $S(t)$ and Remark \ref{r:2.0}.
\begin{proposition}\label{p:1.1} Let $x_n$ be  a bounded sequence of $X$.
Assume that $S(t)x_n$ is contained in a bounded subset $B$ of $X$
for all $t\geq 0$ and $n\geq 1$. For any sequence  $\tau_n\ra+\8$,
define
$$
\gamma_n(t)=S(\tau_n+t)x_n,\Hs t\in(-\tau_n,+\8).
$$
Then there is a subsequence of $\gamma_n$ that converges to a
complete trajectory $\gamma$ of $S(t)$ uniformly on any compact
interval.

\end{proposition}

\subsection{Attractors}

Let $A$ be a  subset of $X$.

 We say that $A$ { attracts} $B\subset
X$, this means that for any $\ve>0$, there exists a $T>0$ such that
$$
S(t)B\subset \cB(A,\ve),\Hs \A\, t>T.
$$
The {attraction basin} of  $A$, denoted by $\W(A)$, is defined as
$$\W(A)=\{x|\,\,\lim_{t\ra\8}d(S(t)x,\,A)=0\}.$$
 The set $A$ is said to be {positively invariant} (resp. {
invariant}), if $$S(t)A\subset A\hs (\,\mb{resp}.
\,\,\,S(t)A=A\,),\Hs \A\,t\geq0.$$
 The { $\omega$-limit set}
$\omega(A)$ of $A$ is defined  as
$$
\omega(A)=\{y\in X|\,\,\,\,\E\, x_n\in A\mbox{ and }
t_n\rightarrow+\infty \mb{ such that } S(t_n)x_n\rightarrow y\}.
$$

Let $\gamma$ be a trajectory on $(t_0,+\8)$ (resp. $(-\8,t_0)$). We
define the $\w$-limit set (resp. $\alpha$-limit set) of $\gamma$ as
follows:
$$
\w(\gamma)=\{y\in X|\,\,\,\,\E\,  t_n\ra+\8 \mb{ such that
}\gamma(t_n)\ra y\}.
$$
$$
\alpha(\gamma)=\{y\in X|\,\,\,\,\E\,  t_n\ra-\8 \mb{ such that
}\gamma(t_n)\ra y\}.
$$
\vs The proofs of  the basic facts listed below can be found in
\cite{Sell}, we omit the details.

\begin{proposition}\label{p:2.2} Let $A\subset X$. If\, $\Cup_{t\geq t_0}S(t)A$  is bounded for some
$t_0>0$, then $\w(A)$ is a nonempty compact invariant set of $S(t)$.
In addition, $\w(A)$ attracts $A$.

If $A$ is connected, then so is $\w(A)$.
\end{proposition}
\begin{definition}
A compact  set $\a\subset X$ is said to be an {attractor} of $S(t)$,
if it is invariant and attracts
 a neighborhood $U$ of itself.

An attractor $\a$ is said to be the {global attractor} of $S(t)$, if
it attracts each bounded subset of $X$.
\end{definition}

\begin{proposition}\label{p:2.0} Let  $\a$ be an
attractor with attraction basin  $\W=\W(A)$. Then
\begin{enumerate}
\item[$(1)$] for any $x\in\W$, there is an $\ve>0$ such that $\a$
attracts $\cB(x,\ve)$, hence $\W$ is open;
\item[$(2)$] $\a$ attracts each
 compact subset $K$ of $\W$;
\item[$(3)$]
$\a$ is Lyapunov stable.
\end{enumerate}
\end{proposition}

\Vs \bt\label{p:2.4} $(${\bf Existence of attractors}\,$)$ Assume
that there is a bounded closed subset $A\subset X$ which attracts a
neighborhood of itself.

Then the semiflow $S(t)$ has an attractor $\a$  contained in $A$.
\et

\subsection{Morse decomposition of attractors}

Let $\a$ be an attractor of $S(t)$. Since $\a$ is invariant, the
restriction $S_\a(t)$ of $S(t)$ on $\a$ is also a semidynamical
system.

A compact set $A$ is said to be an { attractor of $S(t)$ in $\a$},
this  means that $A$ is an attractor of the restricted system
$S_\a(t)$ in $\a$. \bl\label{l:2.9}$\cite{Kap,Lisiam}$ Let  $A$ be
an  attractor of $S(t)$ in $\a$. Then it is  an attractor of $S(t)$
in $X$. \el

For an attractor $A$
 of $S(t)$ in $\a$, define
\be\label{e:3.0a}A^*=\{x\in \a|\,\, \w(x)\cap A=\emptyset\}.\ee
 $A^*$ is said to be the {\bf repeller }of
$S(t)$ in $\a$  dual to $A$, and $(A,A^*)$ is said to be an {\bf
attractor-repeller pair} in $\a$.

A  repeller $A^*$ is invariant; moreover,
 \be\label{e:r1}A^*=\a\setminus \W_\a(A)=\a\setminus \W(A),\ee
where $\W_\a(A)$ is the attraction basin of $A$ in $\a$; see, e.g.,
\cite{Kap,Lisiam}.

\begin{definition}\label{d:3.2}Let $\a$ be an attractor.
An ordered collection $\cM=\{M_1,\cdots,M_l\}$ of subsets of $\a$ is
called a Morse decomposition of $\a$, if there exists an increasing
sequence \be\label{e:3.14} \ep=A_0\varsubsetneq
A_1\varsubsetneq\cdots\varsubsetneq A_l=\a\ee of attractors of $
S(t)$ in $\a$ such that \be\label{e:3.1} M_k=A_k\cap A_{k-1}^*,\Hs
1\leq k\leq l. \ee

\end{definition}

For convenience, the attractor sequence in (\ref{e:3.14}) will be
referred to as a {\bf Morse filtration}, and each $M_k$ in
(\ref{e:3.1})
 a {\bf Morse set}.

 It is
well known that an attractor of a gradient system with finite number
of equilibria has a natural Morse decomposition with each Morse set
being precisely an equilibrium. This was used by many authors to
study continuity of attractors with respect to perturbations
\cite{Cava,CK, ES,Hale, Lidcds,Rob}.

\begin{theorem}\label{t:3.m1}$\cite{Kap}$
Let $\cM=\{M_1,\cdots,M_l\}$ be a Morse decomposition of $\a$ with
Morse filtration $$ \ep=A_0\varsubsetneq
A_1\varsubsetneq\cdots\varsubsetneq A_l=\a.$$ Then
\begin{enumerate}
\item[$(1)$] for each $k$, $(A_{k-1},M_k)$ is an
 attractor-repeller pair in $A_k$;

\item[$(2)$] $M_k$ $(1\leq k\leq l)$ are pair-wise disjoint invariant compact sets;
 \item[$(3)$] if
$\gamma$ is a complete trajectory, then either $\gamma(\R)\subset
M_k$ for some $k$, or else there are indices $i<j$ such that $
\alpha(\gamma)\subset M_j$ and $\w(\gamma)\subset M_i; $
\item[$(4)$] the attractors $A_k$ are uniquely determined by the
Morse sets, that is,
$$
A_k=\Cup_{1\leq i\leq k}W^u(M_i),\Hs 1\leq k\leq l,
$$
where $$ W^u(M_i)=\{\gamma(0)|\,\,\mb{$\gamma$ is a  complete
trajectory in $\a$  with $\alpha(\gamma)\subset M_i\}$}.$$

\end{enumerate}
\end{theorem}

\begin{theorem}\label{t:3.m2} $\cite{Kap}$ Let $\cM=\{M_1,\cdots,M_l\}$ be an ordered collection of
pairwise disjoint compact invariant subsets of $\a$. Suppose that
for every $x\in \a$ and every complete trajectory $\gamma$ through
$x$, we have either $\gamma(\R)\subset M_i$ for some $i$, or else
there are indices $i<j$ such that $\alpha(\gamma)\subset M_j$ and
$\w(\gamma)\subset M_i$.

Then $\cM$ is a Morse decomposition of $\a$.
\end{theorem}

\section{Lyapunov Functions of  Attractors}
In general it is easy to construct  continuous Lyapunov functions of
 attractors for semiflows in metric spaces; see, for instance,
\cite{Kap}. However, to construct global Morse-Lyapunov functions
for Morse decompositions, we need some Lyapunov functions of
attractors which possess nice properties. This is our main concern
in this section.

Let  $\a$ be an attractor, $U$ be a neighborhood of $\a$.
 A
function $V\in C(U)$ is said to be a {\bf Lyapunov function }of $\a$
on $U$, if $V$ is decreasing along trajectories in $U$ with the
restriction $V|_\a$ of $V$ on $\a$ being a constant function.

Let $U$ be an open subset of $X$, and $V\in C(U)$. For $x\in U$,
define
$$
D^+V(x)=D^+_SV(x):=\limsup_{t\ra0^+}\frac{V(S(t)x)-V(x)}{t}.$$
$D^+V(x)$ is called the {\bf Dini derivative} of $V$ along the
semiflow $S(t)$. \vs The main result in this section is contained in
the following theorem.

 \bt\label{t:2.1} $(${\bf Existence of good Lyapunov functions}$)$ Let  $\a$ be an attractor with attraction
basin $\W=\W(\a)$. Then there exists a function $V\in C(X)$
satisfying: \be\label{e:v1} V(x)\equiv 0\,\,(\mb{on }\a), \Hs
V(x)\equiv 1\,\,(\mb{on }X\setminus \W); \ee \be\label{e:v2}
D^+V(x)\leq -v(x),\Hs \A\,x\in X,\ee where
 $v\in C(X)$ is a nonnegative function
satisfying \be\label{e:v3} v(x)>0 \mb{ }(\,x\in \W\setminus \a),\Hs
v(x)= 0 \,\,(\,x\not\in \W\setminus\a). \ee

In case $X=\R^m$, for any $k\in \mathbb{N}$ we can also require
$V\in C^k(X)$.

\et

To construct such a good  Lyapunov function, we first need to do
some auxiliary works.

Let $A,B$ be two subsets of $X$. We will use the notation
``$A\subset\subset B$\,'' to indicate that  $\ol A\subset
\mb{int}\,B$, in addition,
$$\mb{d}(A,\pa B)>0.$$

\begin{definition}   A function
$\alpha\in C(\W)$
  is said to be
{ radially unbounded} on $\W$, if for any $R>0$, there exists a
bounded closed subset $B\ss \W$ such that
$$
\alpha(x)>R,\Hs \A\,x\in \W\-B.
$$

Let $K$ be a bounded closed  subset of $\W$ with $K\ss \W$. \,A
radially unbounded function $\alpha\in C(\W)$ is said to be a {
$\mathcal{K}_0^\8$ function} of $K$ on $\W$, if it satisfies:
$$
\alpha(x)=0\Longleftrightarrow x\in K.$$ A $\mathcal{K}_0^\8$
function $\alpha$ of $K$ is said to be { coercive}, if for any
$\ve>0$,
$$
\alpha(x)\geq\de>0,\Hs \A\, x\in \W\setminus \cB(K,\ve).$$

\end{definition}

\bl\label{l:3.1}Let $\W$ be an open subset of $X$,  $K$ be a bounded
closed subset of $\W$ with $K\ss\W$. Then there exists a coercive
$\mathcal{K}_0^\8$ function $\alpha$ of $K$ on $\W$. \el
\noindent{\bf Proof.} Since $K\ss\W$, we have
$r:=\mb{d}(K,\pa\W)>0$. Set $\W_0=\ol\cB(K,r/2)$. For each $n\in
\mathbb{N}$, define
$$
\W_n=\{x\in\W|\,\,\, d(x,\pa\W)\geq {r}/{(n+2)}\}\Cap \ol\cB(K,nr),
$$
$$
\W_{-n}=\ol\cB\(K,r/2(n+1)\).
$$
Then we obtain a sequence   of bounded closed subsets $\W_k$
\,($k\in \mathbb{Z}$) of $\W$ such that
\begin{enumerate}
\item[(1)] $ \W_{k}\ss \W_{k+1}$ {for all } $k\in \mathbb{Z}$, and
$$
\Cup_{k\in \mathbb{Z}}\W_k=\W,\Hs \Cap_{k\in \mathbb{Z}}\W_k=K;$$
\item[(2)] for any $\ve>0$, we have $\W_{-k}\subset \cB(K,\ve)$ for $k>0$
sufficiently large;  \item[(3)] for any $k\in \mathbb{Z}$, we have
$\cB(K,\ve)\subset \W_k$ for $\ve>0$ sufficiently small.
\end{enumerate}
For each $k\in \mathbb{Z}$ take a continuous function
$\sig_k:\W\ra[0,1]$ such that
$$
\sig_k(x)=0\,\,\,(x\in\W_k),\Hs \sig_k(x)=1\,\,\,(x\in \pa\W_{k+1}).
$$
Pick a sequence $r_k\in (0,+\8)$ ($k\in \mathbb{Z}$) so that
$$
r_k\leq r_{k+1}\,\,\,(\A\, k\in \mathbb{Z}),\Hs \lim_{k\ra-\8}r_k=0.
$$
Let $H_k=\W_k\-\W_{k-1}$ for $k\in \mathbb{Z}$. Since
$\W_{k-1}\subset \mb{int}\W_k$, we have $\pa H_k=\pa\W_{k-1}\cup\pa
\W_k$. Now define $\alpha$ as
$$
\alpha(x)=\left\{\ba{ll}r_{k-1}+\sig_{k-1}(x)(r_k-r_{k-1}),& x\in H_k\mb{ with }k\leq 0;\\[1ex]
(r_{k-1}+k-1)+\sig_{k-1}(x)(r_k-r_{k-1}+1),\hs\hs & x\in H_k\mb{
with
}k>0;\\[1ex]
0,& x\in K.\ea\right.
$$
It is trivial to check that $\a$ is  continuous and is a coercive
$\mathcal{K}_0^\8$ function of $K$ on $\W$. We omit the details of
the argument.

\bl\label{t:2.2} Let $\a$ be an attractor  with attraction basin
$\W=\W(\a)$. Then for any $\de>0$ with $K:=\ol \cB(\a,\de)\ss \W$,
there  exists  a function $V\in C(X)$ such that \be\label{e:v13}
V(x)= 0 \,\,\,(\,x\in\a\,),\Hs 0<V(x)\leq 1\,\,\,(x\not\in \ol
\cB(\a,\de)\,),\ee\be\label{e:v15}\Hs V(x)= 1,\Hs x\in X\setminus
\W;\ee \be\label{e:v14} D^+V(x)\leq -v(x),\Hs x\in X,\ee where
 $v\in C(X)$ is a nonnegative function
satisfying \be\label{e:v4} 0<v(x)\leq 1\,\,\,(\,x\in \W\setminus
\ol\cB(\a,\de)\,),\Hs v(x)=0 \,\,(\,x\not\in\W\setminus \a\,). \ee

If $X=\R^m$, then for any given $k\in \mathbb{N}$, $V$ can be
constructed to be a function of $C^k$.

\el

\noindent{\bf Proof.} Let $\de>0$ be such that $K:=\ol
\cB(\a,\de)\ss \W,$ and let $\alpha$ be a coercive
$\mathcal{K}_0^\8$ function of $K$ on $\W$. Define
$$ \phi(x)=\sup_{t\geq 0}\alpha(S(t)x),\Hs x\in \W.
$$
For each $x\in \W$, since $\a$ attracts a neighborhood
$\cB(x,r)\subset\W$, we have \be\label{e:3.9} S(t)\cB(x,r)\subset
\cB(\a,\de/2),\Hs \A\, t>T \ee for some $T>0$. Thus
$$\phi(y)=\sup_{0\leq t\leq T}\alpha(S(t)y),\Hs y\in \cB(x,r),$$ from which it can be
easily seen that $\phi$ is well defined and continuous on $\W$. We
observe that for any $x\in \W$ and $\tau>0$,
\be\label{e:L1}\phi(S(\tau)x)=\sup_{t\geq
0}\alpha\(S(t)S(\tau)x\)=\sup_{t\geq \tau}\alpha(S(t)x)\leq \phi(x).
\ee By invariance of $\a$ and the definition of $\phi$ we see that
$$\phi(x)=0\,\,\,(x\in \a),\Hs \phi(x)\geq \alpha(x)\,\,\,(x\in \W).$$

Fix a $\lam>0$ and define  $\psi$ on $\W$ as \be\label{e:3.10a}
\psi(x)=\int_0^\8 e^{\lam t}\,\alpha\(S(t)x\)\,dt,\Hs x\in\W. \ee We
deduce by (\ref{e:3.9}) that $\psi$ is well defined and continuous
on $\W$. Again by invariance of $\a$ we have
$$\psi(x)=0,\Hs \A\,x\in \a.$$

Now we evaluate  $D^+V(x)$. Let $\tau>0$. Then
$$
\ba{lll} \psi\(S(\tau)x\)&=\int_0^\8e^{\lam
t}\,\alpha\(S(t)S(\tau)x\)\,dt\\[2ex]
&=e^{-\lam \tau}\int_0^\8e^{\lam
(t+\tau)}\,\alpha\(S(t+\tau)x\)\,dt\\[2ex]
&=e^{-\lam \tau}\int_\tau^\8e^{\lam
t}\,\alpha\(S(t)x\)\,dt\\[2ex]
&\leq e^{-\lam\tau}\psi(x).\ea
$$
Note that this implies \be\label{e:3.13} D^+\psi(x)\leq -\lam
\psi(x). \ee

Set $$L(x)=\phi(x)+\psi(x).$$Then $$D^+L(x)\leq D^+\psi(x)\leq -\lam
\psi(x).$$ Since $ \phi(x)\geq\alpha(x)$, $L$ is radially unbounded
on $\W$. It is clear that
$$
L(x)\equiv 0\,\,\,(\mb{on }\a),\Hs L(x)>0\,\,\,(x\in
\W\setminus\ol\cB(\a,\de)).$$ Let $\mu(s)=1-e^{-s}$ ($s\geq 0$).
Define
$$
V(x)=\left\{\ba{ll}\mu(L(x)), \hs\hs &x\in\W;\\[1ex]
1,&x\in X\setminus \W.\ea\right.$$ We claim that $V$ is continuous
at any point $x_0\in\pa\W$, and hence $V\in C(X)$.

Indeed for any $\ve>0$, there exists an $R>0$ such that
$$1-\ve<\mu(s)<1,\Hs \A\,s\geq R.$$ By radial unboundedness of $L$ one deduces that there exists a bounded closed subset
$K\ss\W$ such that $L(x)>R$ for all $x\in \W\setminus K$. It then
follows that
$$
1-\ve<V(x)=\mu(L(x))<1,\Hs\A\,x\in\W\setminus K,$$ which proves the
claim.

Let $w(x)=d\(x,\,X\-\W\)$. Define
$$
v(x)=\left\{\ba{ll}\min\(\lam\,\mu'(L(x))\psi(x),\,w(x),\,1\),\hs\hs
&x\in
\W;\\[1ex]
0,&x\in X\-\W.\ea\right.$$ Then $v\in C(X)$ and satisfies
(\ref{e:v4}). For $x\in\W$, we observe that
$$
D^+V(x)=\mu'(L(x))D^+L(x)\leq -\lam\,\mu'(L(x))\psi(x)\leq -v(x).
$$
Hence (\ref{e:v14}) holds true.

The proof for the general case is complete. \vs

Now consider the case $X=\R^m$. Pick a sequence of compact subsets
$K_n$  of $\W$ so that
$$
\ol\cB(\a,\de)\ss K_1\ss K_2\ss\cdots,\Hs \W=\Cup_{n\geq1}K_n.
$$
We claim that for each $n\geq 1$, there exists a $\tau_n>0$ such
that \be\label{e:3.11} t_n(x)\geq \tau_n,\Hs \A x\in \pa K_n, \ee
where
$$
t_n(x)=\sup\{t>0|\,\,S([0,t))x\subset \W\setminus K_{n-1}\}. $$
Indeed, if this was not the case, there would exist a sequence
$x_i\in\pa K_n$ such that $s_i:=t_n(x_i)\ra 0$ as $i\ra\8$. By
compactness we can assume that $\lim_{i\ra\8}x_i= x_0$. Of course
$x_0\in \pa K_n$. On the other hand, since $S(s_i)x_i\in K_{n-1}$,
setting $i\ra\8$ one finds
$$x_0=S(0)x_0=\lim_{k\ra\8}S(s_i)x_i\in K_{n-1}.$$ This leads to a
contradiction and proves (\ref{e:3.11}).

Take a nonnegative function $a_0\in C^\8(X)$ with \be\label{e:2.2}
a_0(x)\equiv 0 \mb{ on }\cB(\a,\de/2), \hs \mb{and } a_0(x)>0 \mb{
on }X\setminus \ol\cB(\a,\de).\ee For each $n\geq 2$ we choose a
nonnegative function $a_n\in C^\8(X)$ with
$$
a_n(x)=\left\{\ba{ll}1, \hs\hs&x\in K_n\setminus K_{n-1},\\[1ex]
0, &x\ne K_{n+1}\setminus K_{n-2}.\ea\right.
$$
( $a_n$ can be obtained by appropriately smoothing some continuous
ones.\,) Let \be\label{e:2.1}
\alpha(x)=a_0(x)+\sum_{n=2}^{\8}\frac{n}{\tau_n}a_n(x). \ee Note
that the righthand side of (\ref{e:2.1}) is in fact a finite sum
over $n$. Consequently $\alpha\in C^\8(X)$. Define $\psi(x)$ as in
(\ref{e:3.10a}). Since $\alpha(x)=0$ for $x\in\cB(\a,\de/2)$, by
(\ref{e:3.9}) one easily sees that $\psi$ is well defined; moreover,
$\psi\in C^\8(\W)$ and satisfies (\ref{e:3.13}). It is easy to
verify that
$$ \psi(x)\equiv 0\,\mb{ on }\a,\Hs \psi(x)>0\,\mb{ on
}\W\-\ol\cB(\a,\de).$$

We show that $\psi(x)$ is radially unbounded on $\W$. Let  $x\in
\W\setminus K_n$ ($n\geq 2$), and let
$$
s_n=s_n(x):=\sup\{t>0|\,\,S([0,t))x\subset \W\setminus K_{n-1}\}.
$$
Then
$$
S(t)x\in K_n\setminus K_{n-1},\Hs t\in [s_n-\tau_n,\,s_n),
$$
Therefore
$$  \psi(x)\geq \int_{s_n-\tau_n}^{s_n}\alpha(S(t)x)\,dt\geq \frac{n}{\tau_n}\int_{s_n-\tau_n}^{s_n}a_n(S(t)x)\,dt=n,
$$
and the conclusion follows.

Take a sequence $r_n\ra+\8$ so that
$$
\ol\cB(\a,\de)\subset B_n:=\{x|\,\,\psi(x)\leq r_n\}
$$
for all $n$. Then $ \W=\Cup_{n\geq 1}B_n$. For each $n$, define
\be\label{e:2.5}
\psi_n(x)=\left\{\ba{ll} \psi(x),\hs \hs&x\in B_n;\\[1ex]
r_n,&x\in X\setminus B_n.\ea\right. \ee Clearly $\psi_n\in
C^\8(X\-\pa B_n)$. We will make a slight modification with $\psi_n$
to obtain a smooth function $V_n\in C^k(X)$. For this purpose we
first note that $\psi_n$ is globally Lipschitz on $X$. Therefore
\be\label{s1} |\psi_n(y)-\psi_n(x)|\leq C|x-y|,\Hs
\A\,x\in\pa\W,\,\, y\in X\ee for some $C>0$.
Set $$ G(s)=\left\{\ba{lll}\mb{sgn}\,(s) e^{-1/s^2},\hs\hs& s>0;\\[1ex]
0,& s=0,\ea\right.$$ where sgn($s$) is the signal function. Then
$G\in C^\8(\R)$, and
$$
G^{(n)}(0)=0\,\,\,(\A\, n\geq0),\Hs G'(s)>0\,\,(\,\A\,x\ne0\,).
$$
Define
$$
V_n(x)=G\(\psi_n(x)-r_n\)+G(r_n).
$$
By (\ref{s1}) one easily checks that $V_n\in C^\8(X)$. $V_n$
satisfies:
$$ V_n(x)\equiv 0\,\mb{ on
}\a,\hs V_n(x)>0\,\mb{ on }\W\-\ol\cB(\a,\de).$$ If $x\in B_n$, then
$$ D^+V_n(x)=G'\(\psi_n(x)-r_n\)D^+\psi_n(x)\leq -\lam
G'\(\psi_n(x)-r_n\)\psi(x).
$$
The above estimate naturally holds  for $x\in X\-B_n$, as in this
case both sides equal $0$. \vs

Since $V_n$ is constant on $X\-B_n$, we see that
$\frac{\pa^lV_n(x)}{\pa x_{i_1}\cdots\pa x_{i_l}}$ is bounded on
$X=\R^m$ for any $l\geq0$ and $1\leq i_1,\cdots,i_l\leq m$. Let
$$
c_n=\max_{x\in X}\left|V_n(x)\right|+\sum_{1\leq l\leq k}
\(\sum_{1\leq i_1,\cdots,i_l\leq m} \max_{x\in
X}\left|\frac{\pa^lV_n(x)}{\pa x_{i_1}\cdots\pa x_{i_l}}\right|\).
$$
We may assume that $c_n\geq 1$. Define \be\label{e:3.12}
V(x)=\gamma\sum_{n=1}^\8\frac{1}{2^nc_n}V_n(x),\Hs x\in X, \ee where
$\gamma=1/\sum_{n=1}^\8\frac{1}{2^nc_n}G(r_n)$. Noting that the
series $\sum_{n=1}^\8\frac{1}{2^nc_n}\frac{\pa^lV_n(x)}{\pa
x_{i_1}\cdots\pa x_{i_l}}$ is uniformly convergent on $X$ for any
$0\leq l\leq k$ and $1\leq i_1,\cdots,i_l\leq m$, one concludes that
$V\in C^k(X)$.

It is trivial to check that $V$ satisfies all the other properties
required in the lemma.

The proof of the lemma is complete.

\Vs\noindent {\bf Proof of Theorem \ref{t:2.1}. } Take a sequence of
positive numbers
$$\de_0>\de_1>\cdots >\de_n\ra 0$$ with
$\ol\cB(\a,\de_0)\ss \W$. Then for each $\de_n$  one can find
functions $V_n,v_n\in C(X)$ satisfying all the properties in Lemma
\ref{t:2.2} with $V,v$ and $\de$ therein replaced by $V_n,v_n$ and
$\de_n$, respectively. Define
$$
V(x)=\sum_{n=0}^\8\frac{1}{2^{n+1}}\,V_n(x),\Hs x\in X.
$$
 Then $V$ is a function satisfying (1) and (2) in Theorem \ref{t:2.1} with
$$v(x)=\sum_{n=0}^\8\frac{1}{2^{n+1}}\,v_n(x).$$

In case $X=\R^m$, one can define a function $V\in C^k(\W)$ in the
same manner as in (\ref{e:3.12}). We omit the details.

\vs\Vs

The following result can be obtained directly from Theorem
\ref{t:2.1}. It is  also readily  implied in the  proof of Lemma
\ref{t:2.2}.

 \bt\label{t:2.0} Let $\a$ be an attractor
with attraction basin $\W=\W(\a)$. Then $\a$ has a radially
unbounded Lyapunov function $L$ on $\W$ such that
\be\label{e:v9}L(x)\equiv 0\,\,(\mb{on }\a\,), \Hs D^+L(x)\leq
-v(x)\,\,(\A\,x\in \W),\ee
 where
 $v\in C(\W)$ is a nonnegative function
satisfying \be\label{e:v3} v(x)>0 \mb{ }(\mb{for }\,x\in \W\setminus
\a),\Hs v(x)\equiv0 \,\,(\mb{on }\a). \ee

If $X=\R^m$, then for any given $k\in \mathbb{N}$, $L$ can be
constructed in $C^k(\W)$.

\et

\noindent{\bf Proof.} Let $V$ be a Lyapunov function of $\a$ given
by Theorem \ref{t:2.1} Define $$L(x)=\eta(V(x)),\Hs x\in\W,$$ where
 $\eta(s)=-\ln(1-s)$ ($s\in[0,1)$). Then $L$ is a radially unbounded
Lyapunov function on $\W$ that satisfies all the desired properties.
\Vs The validity of the conclusions in the following remarks can be
easily seen from the proofs of Theorems \ref{t:2.1} and \ref{t:2.0}.

\br\label{r:2.1}If we replace the attraction basin $\W$ in Theorems
\ref{t:2.1} and \ref{t:2.0} by any positively invariant open
neighborhood $U$ of $\a$, then all the conclusions in the theorems
remain valid.\er

\br\label{r:2.2} If we replace the attractor $\a$ by any positively
invariant bounded closed neighborhood $W$ of $\a$ with $W\ss\W$,
then all the conclusions in Theorems \ref{t:2.1} and \ref{t:2.0}
hold true.\er

\section{Morse-Lyapunov Functions}

Let there be given an attractor $\a$ with attraction basin
$\W=\W(\a)$ and  Morse decomposition $\cM=\{M_1,\cdots,M_l\}$. Let
$
\cR=\Cup_{1\leq k\leq l}M_k.
$
\begin{definition} $V\in C(\W)$ is said to be a $($global$\,)$ {Morse-Lyapunov function}
$($M-L function in short\,$)$ of $\cM$ on $\W$,
if\begin{enumerate}\item[$(1)$] $V$ is constant on each Morse set
$M_k$; \item[$(2)$] $V(S(t)x)$ is strictly decreasing in $t$ for
$x\in \W\-\cR$.\end{enumerate} $V$ is said to be a strict M-L
function of $\cM$, if in addition it satisfies
$$V(M_1)<V(M_2)<\cdots<V(M_l).$$
\end{definition}

\br If we restrict the system on the attractor $\a$, then  a M-L
function of $\cM$ can be easily formulated as in \cite{Conley}. See
also \cite{Mis}. \er The main  result in this section is contained
in the following theorem.

\bt\label{t:3.1} $(${\bf Existence of strict M-L functions}$)$ $\cM$
has a radially  unbounded strict M-L function $V$ satisfying
\be\label{e:v8} D^+V(x)\leq -v(x),\Hs\mb{for } x\in \W,\ee where
 $v\in C(\W)$ is a nonnegative function
satisfying \be\label{e:3.7} v(x)\equiv 0 \,\,\,(x\in
\mathcal{R}\,),\Hs v(x)>0\,\,\,( x\in \W\setminus \mathcal{R}\,).
\ee

In case $X=\R^m$, for any given $n\geq 1$ we can also construct $V$
so that  $V\in C^n(\W)$.

\et

\noindent{\bf Proof.} Let $ \ep=A_0\varsubsetneq
A_1\varsubsetneq\cdots\varsubsetneq A_l=\a$ be the Morse filtration
of $\cM$.

For $k=l$, we infer from Theorem \ref{t:2.0} that there is   a
 radially unbounded nonnegative function $V_l\in C(\W)$ such that
\be V_l(x)=0,\Hs \A\,x\in A_l=\a,\ee \be\label{e:3.3}D^+V_l(x)\leq
-v_l(x),\Hs\A\, x\in \W,\ee where
 $v_l\in C(\W)$ is a nonnegative function
satisfying (\ref{e:v3}).

For each $1\leq k\leq l-1$,  by Theorem \ref{t:2.1}  there exists a
$V_k\in C(X)$ such that \be\label{e:3.4} V_k(x)\equiv 0\,\mb{ (on
}A_k),\Hs V_k(x)\equiv 1\,\mb{ (on }X\setminus \W(A_k));\ee
\be\label{e:3.6}D^+V_k(x)\leq -v_k(x),\Hs \A\,x\in X,\ee where
 $v_k\in C(X)$ is a nonnegative function
with \be\label{e:v12} v_k(x)>0 \,\,\mb{ (}x\in \W(A_k)\setminus
A_k\,),\,\Hs v_k(x)=0 \mb{ (}x\not\in \W(A_k)\-A_k\,). \ee Define
$V\in C(\W)$ as: \be\label{e:5.18} V(x)=\sum_{1\leq k\leq l}V_k(x),
\Hs x\in\W. \ee We show that $V$ has all the required properties
with $v(x)=\sum_{1\leq k\leq l}v_k(x).$

First we recall that
$$
M_k=A_k\cap
A_{k-1}^*=A_k\cap\(\a\setminus\W(A_{k-1})\)=A_k\cap\W(A_{k-1})^c,
$$
where $\W(A_{k-1})^c=X\setminus \W(A_{k-1})$. Thus if $i\leq k-1$,
then $ M_k\subset  \W(A_{k-1})^c\subset \W(A_i)^c. $ It follows that
$$V_i(M_k)=1,\Hs\mb{for } 1\leq i\leq k-1.$$ On the other hand if $i\geq k$, then
$M_k\subset A_k\subset A_i$. Therefore we find
$$V_i(M_k)=0,\Hs\mb{for }  i\geq k.$$ Hence we deduce that
$$
V(M_k)=\sum_{1\leq i\leq l}V_i(M_k)=k-1.
$$

We  observe that
$$
D^+V(x)\leq \sum_{1\leq i\leq l}D^+V_i(x)\leq -\sum_{1\leq i\leq
l}v_i(x)=-v(x),\Hs x\in\W.
$$

There remains to check  (\ref{e:3.7}).

Let $x\in\cR$. We may assume that $x\in M_k=A_k\cap A_{k-1}^*$. If
$i\geq k$, then we have $x\in A_i$ for all $i\geq k$, therefore
$v_i(x)=0$. If $i\leq k-1$, then we deduce by $x\in A_{k-1}^*$ that
$x\not\in \W(A_{k-1})\supset \W(A_i)$, which  implies $v_i(x)=0$. In
conclusion, $v_i(x)=0$ for all $1\leq i\leq l$. Consequently
$v(x)=0$.

Now assume $x\in \W\setminus \mathcal{R}$. Then there is a smallest
$k$ such that $x\in \W(A_k)$. We claim that $x\not\in A_k$. Indeed,
if $x\in A_{k}$, then either $x\in M_k$, or $x\in\W(A_{k-1})$. In
any case one will get a contradiction. Hence the claim holds true.
Since $x\in \W(A_k)\-A_k$, we see that $v_k(x)>0$, and hence
$v(x)\geq v_k(x)>0.$

In case $X=\R^m$, each function $V_k$ can be constructed in $C^n$.
Consequently $V\in C^n(\W)$. The proof is complete. \Vs

The following  interesting result indicates that we can modify any
M-L function to a strict one by ``preserving\,'' its  values in
small neighborhoods of Morse sets. Therefore for many purposes it
suffices to consider strict M-L functions.

\bt\label{t:3.2} There exists an $\ve>0$ sufficiently small such
that for any M-L function $V$ of $\cM$, one can find a strict M-L
function $L$ of $\cM$ such that \be\label{e:4.21}(L-V)(x)\equiv
const.\ee on the $\ve$-neighborhood $\cB(M_k,\ve)$ of each Morse set
$M_k$.

\et

\noindent{\bf Proof.} Let $M_k$ be any given Morse set. To prove the
result, it suffices to show that there exist an $\ve>0$
 and an M-L function $L$ of $\cM$ such that (\ref{e:4.21}) holds
 true, moreover,
\be\label{e:4.22} L(M_{j})>L(M_k)>L(M_{i}),\Hs \A\,j>k>i. \ee

Take a positively invariant bounded closed neighborhood $W$ of $A_k$
with $W\ss\W(A_k)$ (this can be done by using  Lyapunov functions of
$A_k$). By Remark \ref{r:2.2} and Theorem \ref{t:2.0} we deduce that
there exists a $\mathcal{K}_0^\8$ function $\Phi$ of $W$
 on $\W(A_{k})$ such that  $$ \Phi(x)\equiv
0\,\,\,(x\in W), \Hs D^+\Phi(x)\leq 0\,\,\,(x\in \W(A_{k})\,).$$ Fix
an $R>0$ so that
$$
W\subset \Phi_R:=\{x\in \W(A_{k})|\,\,\Phi(x)\leq R\}\ss \W(A_{k}).
$$
Let
$$
V_1(x)=\min\(\frac{1}{R}\Phi(x),\,1\),\Hs x\in X.
$$
Then $V_1$ satisfies \be\label{e:4.17} V_1(x)=0\,\,(x\in W),\Hs
V_1(x)=1\,\,\,(x\in X\-\Phi_R). \ee Clearly $D^+V_1(x)\leq 0$ for
all $x\in X$.

Since $\Phi_R\ss \W(A_{k})$, we deduce that
$$\mb{d}(M_i,\Phi_R)>0,\Hs \mb{for all }i>k.$$ Therefore if $\ve>0$ is taken
sufficiently small then one has
$$\cB(M_i,\ve)\subset  W\,\,\,(\mb{for } i\leq k),\Hs \cB(M_i,\ve)\subset X\-\Phi_R\,\,\,(\mb{for } i> k). $$
It follows by (\ref{e:4.17}) that \be\label{e:4.18}
V_1\left|_{\mb{\footnotesize B}(M_i,\ve)}\right.\equiv
0\,\,\,(\mb{for } i\leq k),\Hs V_1\left|_{\mb{\footnotesize
B}(M_i,\ve)}\right.\equiv 1\,\,\,(\mb{for } i>k). \ee

Now we choose a positively invariant closed neighborhood $U$ of
$A_{k-1}$ so that
$$U\ss\W(A_{k-1})\cap W.$$ Again by Remark \ref{r:2.2} and Theorem \ref{t:2.0} there exists a $\mathcal{K}_0^\8$ function $\Psi$ of
$U$ on  $\W(A_{k-1})$ such that \be\label{e:3.10} \Psi(x)\equiv
0\,\,(x\in U), \Hs D^+\Psi(x)\leq 0\,\,(x\in \W(A_{k-1})\,).\ee Take
a positive number $a>0$ so that
$$
U\subset \Psi_a:=\{x\in \W(A_{k-1})|\,\,\Psi(x)\leq a\}\ss
\W(A_{k-1}).
$$
Let
$$
V_2(x)=\min\(\frac{1}{a}\Psi(x),\,1\),\Hs x\in X.
$$
Then $V_2$ satisfies \be\label{e:4.19} V_2(x)=0\,\,(x\in U),\Hs
V_2(x)=1\,\,(x\in X\-\Psi_a). \ee Moreover, we have  $D^+V_2(x)\leq
0$ for all $x\in X$.

We further restrict $\ve>0$ sufficiently small so that $$
\cB(M_i,\ve)\subset U\subset W\,\,\,(\mb{for } i<k),\Hs
\cB(M_i,\ve)\subset X\-\Psi_a\,\,\,(\mb{for } i\geq k). $$ Then
(\ref{e:4.19}) implies that \be\label{e:4.20}
V_2\left|_{\mb{\footnotesize B}(M_i,\ve)}\right.\equiv
0\,\,\,(\mb{for } i<k),\Hs V_2\left|_{\mb{\footnotesize
B}(M_i,\ve)}\right.\equiv 1\,\,\,(\mb{for } i\geq k).\ee

Now for any M-L function $V$, define
$$L(x)=V(x)+V_1(x)+V_2(x),\Hs x\in X.
$$
Then by (\ref{e:4.18}) and (\ref{e:4.20})  one finds that for $x\in
\cB(M_i,\ve)$,
$$
L(x)=\left\{\ba{ll}V(x)+2,\hs\hs &\mb{if }
i>k;\\[1ex]V(x)+1,\Hs &\mb{if }
i=k;\\[1ex]V(x),& \mb{if }i<k.\ea\right.
$$
Since $D^+L(x)\leq D^+V(x)$, we conclude immediately that $L$ is a
M-L function and satisfies (\ref{e:4.22}).

The proof of the theorem is complete.
\section{Deformation Lemmas}

In this section we extend some well known  deformation results in
the critical-point theory to semiflows on metric spaces.

Let $E$ be a topological space. A subset   $F\subset E$ is said to
be a { strong deformation retract} of $E$, if there exists a
continuous mapping  $H:[0,1]\X E\ra E$ such that
$$
H(0,\.)|_E=\mb{id}_E,\Hs H(1,E)\subset F,
$$
$$
H(\sig,\.)|_F=\mb{id}_F,\Hs \A\,\sig\in[0,1].
$$

Let $\a$ be an attractor of the system $S(t)$
 with attraction basin $\W=\W(\a)$ and Morse decomposition $\cM=\{M_1,\cdots,
 M_l\}$, and let $V$ be a strict M-L function of $\cM$. For convenience in statement we will call
 $$c_k=V(M_k),\Hs k=1,2,\cdots,l$$ the {\bf generalized critical
 values} of $V$.

 For $a\in\R$ we denote by $V_a$ the level set of $V$ in $\W$,
 $
 V_a=\{x\in\W|\,\,V(x)\leq a\}.
 $
 $V_a$ is clearly positively
 invariant, in addition, one easily checks that $$V_a\subset
 \W(A_k),\Hs\mb{if }a<c_{k+1}.$$

 \bt\label{t:4.1} $(${\bf First Deformation Lemma}$)$ Let $-\8<a<b<+\8$. If $V$ has no generalized critical
 values in $[a,b]$, then $V_a$ is a strong
 deformation retract of $V_b$.\et

\noindent{\bf Proof.} Note that  $V$ attains its minimum on $M_1$.
Therefore  $V_a=V_b=\ep$ if $a,b<c_1$. So it can be  assumed that
$c_k<a<b<c_{k+1}$ for some $1\leq k\leq l$ (set $c_{l+1}=+\8$ for
convenience). Recall that $V_b\subset \W(A_{k})$.

Define a function $t(x)$ on $V_b$ as \be\label{e:4.1}
t(x)=\left\{\ba{ll}\sup\{t\geq0|\,\,S([0,t))x\subset V_b\setminus
V_a\},\hs\hs &x\in V_b\setminus V_a;\\[1ex]
0,& x\in V_a.\ea\right. \ee We first show that $t(x)<+\8$  for any
$x$. Suppose that $t(x)=+\8$ for some $x\in V_b$. Then $V(S(t)x)>a$
for all $t>0$. It follows that $V(y)\geq a$ for all $y\in \w(x)$. On
the other hand since $V_b\subset \W(A_k)$, we find that
$\w(x)\subset A_{k}$, therefore $$V(y)\leq V(M_k)=c_k<a,\Hs\A\,
y\in\w(x),$$ which leads to a contradiction!

Note that  $S(t(x))x\in V_a$. \bl\label{l:4.2}$t(x)$ is a continuous
function of $x$ on $V_b$.\el

\noindent{\bf Proof.} We first consider $x_0\in V_b\setminus V_a$.
Assume $x_n\ra x_0$ as $n\ra \8$. Then $x_n\in V_b\setminus V_a$ for
$n$ sufficiently large. Let $\ve>0$ be given arbitrary. We prove
that for some $n_0>0$,
$$
t(x_n)>t(x_0)-\ve,\Hs \A\,n>n_0.
$$
 Suppose the contrary. Then there exists
a subsequence of $x_n$ (still denoted by $x_n$) such that
$t(x_n)\leq t(x_0)-\ve$ for all $n$. We may assume $t(x_n)\ra
t_0\leq t(x_0)-\ve.$ Letting $n\ra\8$ one finds that
$$
S(t_0)x_0=\lim_{n\ra \8}S(t(x_n))x_n\in V_a,
$$
which leads to a contradiction.

In what follows  we show that
$$
t(x_n)\leq t(x_0)+\ve
$$
for sufficiently large $n$. Indeed, if this was not the case, there
would exists a subsequence ${n_k}$ such that $t(x_{n_k})\geq
t(x_0)+\ve$ for all $k$. Passing to the limit one finds that
 $$S\([0,\,t(x_0)+\ve]\)x_0\subset V^{-1}([a,b]).$$ However, since
$S(t(x_0))x_0\in V_a$ and $S(t(x_0))x_0\not\in\cR$, we see that
$$
V\(S(t(x_0)+\de)x_0\)<V(S(t(x_0))x_0)\leq a
$$
for $\de>0$ sufficiently small, which yields a contradiction.

Now we consider the case $x_0\in V_a$ in which  we have $t(x_0)=0$.
Let $x_n\ra x_0$. We need to prove that $t(x_n)\ra0$. If $x_n\in
V_a$, then by definition of $t(x)$ one has $t(x_n)=0$. So it can be
assumed that $x_n\in V_b\setminus V_a$. Consequently we have
$V(x_0)=a$. If $\lim_{n\ra\8}t(x_n)\ne 0$, then there is a
subsequence of $t(x_n)$, still denoted by $t(x_n)$, such that
$t(x_n)\geq\tau>0$ for all $n$. Since $S([0,\tau])x_n\subset
V_b\setminus V_a$ and $S(t)x_n$ converges to $S(t)x_0$, one deduces
that $S([0,\tau])x_0\subset V_b\setminus V_a$, which implies that
$$a=V(x_0)\geq V(S(t)x_0)\geq a,\Hs t\in[0,\tau],$$ that is, $V(S(t)x_0)\equiv
a$ on $[0,\tau]$. On the other hand, since $x_0\not\in \cR$, we
should have  $$V(S(t))x_0<V(x_0)=a$$ for $t>0$ sufficiently small.
This is a contradiction which finishes the proof of the lemma.

\Vs Now we continue to prove Theorem \ref{t:4.1}. Define
$$
H(\sig,x)= S(\sig\, t(x))x,\Hs x\in V_b.$$ Then $H: [0,1]\X V_b\ra
V_b$ satisfies: $$H(0,\.)=\mb{id}_{V_b},\Hs H(1,V_b)\subset V_a,$$
$$
H(\sig, x)=x,\Hs \A \,(\sig,x)\in [0,1]\X V_a.
$$
Since $t(x)$ is continuous on $V_b$, we see that $H$ is a continuous
mapping.

The proof of Theorem \ref{t:4.1} is complete. \Vs

It is interesting to consider the case $M_k$ is an equilibrium of
the flow. In this case we have the following stronger conclusion.

\bt\label{t:4.3} $(${\bf Second Deformation Lemma}$)$ Assume that
for some $1\leq k_0\leq l$ the  Morse set $M_{k_0}$ consists of
exactly one equilibrium $E$ of the flow. Let $c=c_{k_0}:=V(E)$, and
let $b>c$ be a number such that $V$ has no generalized critical
values in $(c,b]$.

Then $V_c$ is a strong
 deformation retract of $V_b$.\et

\noindent{\bf Proof.} Define a function $t(\sig,x)$ on $[0,1)\X V_b$
as \be\label{e:4.2}
t(\sig,x)=\left\{\ba{ll}\sup\{t\geq0|\,\,S([0,t])x\subset
V_b\setminus
V_{\sig c+(1-\sig)b}\},\hs\hs&x\in V_b\setminus V_{\sig c+(1-\sig)b}\,;\\[1ex]
0,& x\in V_{\sig c+(1-\sig)b}.\ea\right. \ee Then $t(\sig,x)<+\8$
for any $(\sig,x)\in [0,1)\X V_b$, and $S(t(\sig,x))x\in V_{\sig
c+(1-\sig)b}$. By slightly modifying the proof of Lemma \ref{l:4.2}
it can be  easily shown that $t(\sig,x)$ is continuous on $[0,1)\X
V_b$. We observe that $t(\mu,x)$ is nondecreasing in $\mu$.
Therefore the limit
\be\label{e:4.6}\lim_{\mu\ra1^-}t(\mu,x)=T(x)\ee exists.

Define $H:[0,1]\X V_b\ra V_b$ as
$$
H(\sig,x)=\left\{\ba{ll}S(\sig\, t(\sig,x))x,\Hs &\sig<1,\,\,x\in V_b;\\[1ex]
\lim_{\mu\ra 1^-}S(\mu \,t(\mu,x))x,\hs\hs&\sig=1,\,\,x\in V_b\,
.\ea\right.$$ We first show that $$\lim_{\mu\ra 1^-}S(\mu
\,t(\mu,x))x=\lim_{t\ra T(x)}S(t)x\in V_c$$ do exist, hence $H$ is
well defined.

If $T(x)<+\8$, then clearly $\lim_{t\ra T(x)}S(t)x=S(T(x))x.$ So we
assume  $T(x)=+\8$. In this case we necessarily have
$S([0,+\8))x\subset V_b\setminus V_c.$ Noting  that for any
$\sig<1$, $$S(t)x\in V_{\sig c+(1-\sig)b},\Hs \A\,t>t(\sig,x),$$ one
easily deduces that $\lim_{t\ra+\8}V(S(t)x)=c$, and consequently
$V(\w(x))=c$. Because $\w(x)\subset \cR:=\Cup_{1\leq i\leq l}M_i$,
we conclude that $\w(x)=M_{k_0}=E$, that is,
\be\label{e:4.7}\lim_{t\ra T(x)}S(t)x=E.\ee

It is trivial to examine that   $H$ satisfies:
$$H(0,\.)=\mb{id}_{V_b},\Hs H(1,V_b)\subset V_c,$$
$$
H(\sig, x)=x,\Hs \A \,(\sig,x)\in [0,1]\X V_c.
$$
To complete the proof of the theorem, there remains to verify the
continuity of $H$.

By the definition of $H$ it is clear that $H$ is continuous on
$[0,1)\X V_b$, so we only consider the continuity of $H$ at any
point $(1,x_0)\in \{1\}\X V_b$, at which we have
$$H(1,x_0)=\lim_{t\ra T(x_0)}S(t)x_0 .$$

 \vs {\sl Case }1) \,\,``\,$H(1,x_0)=E,\,\,x_0\not\in
 V_c$\,''.
\vs
 This is the worst case we meet, in which one necessarily has $T(x_0)=+\8$.  Let $(\sig_n,x_n)\ra
(1,x_0)$. Then $x_n\in V_b\setminus V_c$ for $n$ sufficiently large.
By some similar argument as in the proof of Lemma \ref{l:4.2} it can
be shown that $T(x_n)\ra +\8$. Consequently
\be\label{e:ts}s_n:=\sig_n\,t(\sig_n,x_n)\ra+\8,\Hs \mb{as
}n\ra\8.\ee (Here we set $t(\sig_n,x_n)=T(x_n)$ if $\sig_n=1$.) We
need to prove that \be\label{e:4.9} H(\sig_n,x_n)=S(s_n)x_n:=y_n\ra
E. \ee

Let $\ve>0$ be given so that $\cB(M_i,\ve)\Cap \cB(M_j,\ve))=\ep$
for all $i\ne j$. Let $r>0$ be such that $\a$ attracts
$B:=\ol\cB(\a,r)\subset\W$. Then for some  $t_0>0$,
$$S(t)B\subset B,\Hs t\geq t_0.$$
Let $N=\ol{\Cup_{t\geq t_0}S(t)B}.$ Take a $\de_0>0$ sufficiently
small so that $\a$ attracts $\cB(x_0,\de_0)$. Then there exists a
$\tau>0$ such that $S(t)\cB(x_0,\de_0)\subset N$ for $t\geq \tau$.
Since $x_n\ra x_0$,
 we can assume that $x_n\in\cB(x_0,\de_0)$ for all $n\geq 1$, and
 hence
 $$
 S(t)x_n\in N,\Hs t\geq \tau,\,\,\,n\in \mathbb{N}.
 $$

We claim that there exists a sequence $t_n$ with $s_n/2<t_n<s_n$
such that \be\label{e:4.13}S(t_n)x_n\ra E\ee as $n\ra\8$. Indeed, if
this was not the case,  one would find  a subsequence  $n_k$  and an
$\eta>0$ such that \be\label{e:4.10}d\(S(t)x_{n_k},\,E\)\geq\eta,\Hs
t\in(s_n/2,s_n).\ee It can be assumed  that $s_n/2> \tau$ for all
$n$. Define a sequence of trajectories $\gamma_k$ as:
$$
\gamma_{k}(t)=S\(t+3s_{n_k}/4\)x_{n_k},\Hs -s_{n_k}/4<t<s_{n_k}/4.$$
Then $\gamma_k$ is contained in $N$. By virtue of Proposition
\ref{p:1.1} one easily deduces that there is a subsequence of
$\gamma_k$ which converges uniformly on any compact interval to a
complete trajectory $\gamma$ of the flow contained in $N$.
(\ref{e:4.10}) then implies that \be\label{e:4.11}
d(\gamma(t),E)\geq\eta>0,\Hs \A\, t\in\R.\ee It is also clear that
$$c\leq V(\gamma(t))\leq b,\Hs t\in\R.$$

On the other hand, since $\gamma$ lies in the attractor $\a$ and
$V_b\subset \W(A_{k_0})$, we deduce that $\gamma$ is contained in
$A_{k_0}$. It then follows by $V(\gamma(t))\geq c=V(E)$ that
$\gamma(t)\equiv E$. This contradicts (\ref{e:4.11}) and proves
(\ref{e:4.13}).

Now we show that \be\label{e:4.14}
d\(H(\sig_n,x_n),\,E\)=d\(S(s_n)x_n,\,E\)\leq \ve\ee for $n$ large
enough, which completes the proof of the continuity of $H$ at
$(1,x_0)$.

Suppose the contrary. Then there exists a subsequence of $n$, which
we relabel as $n$, such that
$$d\(S(s_{n})x_{n},\,E\)>\ve$$ for all $n$. By
(\ref{e:4.13}) we deduce that there exists an $n_0>0$ such that for
each $n>n_0$, one can find an interval $[\theta_n,\,\tau_n]\subset
[t_{n},\,s_{n}]$ such that \be\label{e:4.16}
d\(S(\theta_n)x_{n},\,E\)=\ve/2,\Hs d\(S(\tau_n)x_{n},\,E\)=\ve,\ee
and that $$S(t)x_{n}\in N\cap
\(\ol\cB(E,\ve)\setminus\cB(E,\ve/2)\),\hs \mb{for }t\in
(\theta_n,\tau_n).$$ We claim that there exist $T,\sig>0$ such that
\be\label{e:4.15}0<\sig\leq\tau_n-\theta_n\leq T<+\8\ee for all $n$.
Indeed, if $\tau_n-\theta_n$ is unbounded from above, then  by some
similar argument as above one will find a complete trajectory
$\gamma$ contained in $\ol\cB(E,\ve)\setminus\cB(E,\ve/2)$ with
(\ref{e:4.11}) holds, which leads to a contradiction. Now suppose
that there is a subsequence $n_k$ such that $
\tau_{n_k}-\theta_{n_k}\ra 0$ as $n_k\ra\8$. Since $\theta_n\geq
t_{n}$ and $t_n\ra+\8$, by asymptotic compactness of $S(t)$ it can
be assumed that $S(\theta_{n_k})x_{n_k}\ra y$. Then
$$
S(\tau_{n_k})x_{n_k}=S(\tau_{n_k}-\theta_{n_k})S(\theta_{n_k})x_{n_k}\ra
y,\Hs \mb{as } n_k\ra\8,$$ that is,
$\lim_{k\ra\8}S(\theta_{n_k})x_{n_k}=\lim_{k\ra\8}S(\tau_{n_k})x_{n_k}=y$.
Passing to the limit in (\ref{e:4.16}) for the subsequence $n_k$, it
yields
$$
d\(y,\,E\)=\ve/2,\Hs d\(y,\,E\)=\ve,
$$
a contradiction! Hence (\ref{e:4.15}) holds true.

We now prove that there exists an $r>0$ such that \be\label{e:4.12}
V(S(\tau_n)x_{n})\leq V\(S(\theta_n)x_{n}\)-r,\Hs n=1,2\cdots. \ee
Indeed, if (\ref{e:4.12}) fails to be true, there would exist a
subsequence of ${n}$ (still denoted by ${n}$) such that
$$\lim_{k\ra\8}\left[V\(S(\theta_n)x_{n}\)-V(S(\tau_n)x_{n})\right]= 0.$$
We may assume that $\tau_n-\theta_n\ra \sig>0$, and that
$S(\theta_n)x_{n}\ra y$ (recall that $\theta_n\ra+\8$). Passing to
the limit in the above equation one obtains $V(S(\sig)y)=V(y)$,
which leads to a contradiction and proves (\ref{e:4.12}).

Now by  (\ref{e:4.12}) and the choice of $\theta_n$ we have
$$
V(S(\tau_n)x_{n})\leq V\(S(\theta_n)x_{n}\)-r\leq
V\(S(t_{n})x_{n}\)-r.
$$It then follows by (\ref{e:4.13}) that
$$V(S(\tau_n)x_{n})<c-r/2$$ for $n$ sufficiently large. This
contradicts to the fact that $ V(S(t)x_n)\geq c$ for $t\leq s_n$ and
completes the proof of (\ref{e:4.14}).

\Vs {\em Case} 2)\,\,``\,$H(1,x_0)=E,\,\,x_0\in
 V_c$\,''.
\vs
 In this case by
definition of $H$ we must have $x_0=E$. Let $(\sig_n,x_n)\ra
(1,x_0)$. Note that if $x_{n_k}\in V_c$, then $x_{n_k}\ra x_0$
implies $H(\sig_{n_k},x_{n_k})=x_{n_k}\ra x_0$. Thus it can be
assumed that $x_n\in V_b\setminus V_c$.

 If there exists a $T>0$ such that $T(x_n)\leq T<+\8$ for
all $n$, then the sequence $s_n$ in (\ref{e:ts}) is bounded.
Therefore we directly have
$$
\lim_{n\ra\8}H(\sig_n,x_n)=\lim_{n\ra\8}S\(s_n\)x_n=\lim_{n\ra\8}S\(s_n\)x_0=E.
$$

If $T(x_n)\ra+\8$, then we come back to a situation as in Case 1).

The general case in which  neither $T(x_n)$ is bounded nor
$T(x_n)\ra+\8$ can be treated by a simple contradiction argument.

\Vs {\em Case} 3)\,\,``\,$H(1,x_0)\ne E$\,''. \vs

This is the simplest case in which $T(x_0)<+\8$. Since the argument
is an easy excise, we omit the details. The proof of the theorem is
finished.

\vskip12pt

Now we define an equivalence relation ``$\sim$'' on $X$ as follows:
 $$\mb{If $x\ne y$, then $x\sim y$ $\Longleftrightarrow$ $x$ and $y$ belong to the same
Morse set.}$$ Denote by $[x]$ the equivalence class of $x$, and let
$\stackrel{\sim}{X}=X/\sim$ be the quotient space. Then ``$\sim$''
collapses each Morse set $M_k$ to one point of $\~X$, simply denoted
by $[M_k]$.

Let $q:X\ra\~X$ be the quotient map, i.e., $q(x)=[x]$ for each $x\in
X$. Set
$$
\~S(t)[x]=q\circ S(t)x,\Hs x\in\~X,\,\,t\geq 0.
$$
Then $\~S(t)$ is a semiflow on $\~X$, which will be referred to as
the {quotient semiflow}. $\~S(t)$ has a corresponding attractor
$\~\a=q(\a)$ and Morse decomposition
$$
\~\cM=\left\{[M_1],\,\cdots,[M_l]\right\}.
$$
 Define $\~V([x])=V(x)$ (for $[x]\in \~\W:=q(\W)$).
Since $V$ is constant on each Morse set $M_k$, the function $\~V$ is
well defined and continuous on $\~\W$. It is easy to see that $\~V$
is a strict M-L function of $\~\cM$.

\vs As a direct consequence of Theorem \ref{t:4.3}, we have
\bt\label{t:4.4}  Let $c=c_k$ be a generalized critical value of
$V$, and let $b>c$ be a number such that $V$ has no generalized
critical value in $(c,b]$. Then $\~V_c$ is a strong
 deformation retract of $\~V_b$.\et

\section{Morse Theory of Attractors}
\subsection{Critical groups of Morse sets}

In this subsection we  introduce and discuss the concept of critical
groups for Morse sets.

We will denote by $H_*$ the usual singular homology theory with
coefficients in a given Abelian group $\mathscr{G}$. Let $\a$ be an
attractor of the system $S(t)$
 with attraction basin $\W=\W(\a)$, and let  $\cM=\{M_1,\cdots,
 M_l\}$ be a Morse decomposition of $\a$ with Morse filtration $$\ep=A_0\nsubseteq
A_1\nsubseteq\cdots\nsubseteq A_l=\a,$$  $V$ be a strict M-L
function of $\cM$. Set
$$
c_k=V(M_k),\Hs 1\leq k\leq l.
$$
\begin{definition}\label{d:5.1}Take two numbers $a<b$ be such that  $c_k$ is the unique generalized critical value of $V$ in $[a,b]$.
Then the {\bf critical group} $C_*(M_k)$ of the Morse set $M_k$ is
defined to be the homology theory  given by
$$
C_q(M_k)=H_q\(V_b,\,V_a\),\Hs q=0,1,\cdots. $$
\end{definition}

By virtue of the First Deformation Lemma, one easily understands
that the definition of the critical group does not depend on the
choice of the numbers $a$ and $b$. In case $M_k$ consists of exactly
an equilibrium $E$ of the semiflow,  we even have the following
stronger conclusion.

\begin{proposition}\label{t:5.1} Assume that $M_k$ consists of exactly  an equilibrium $E$ of the
semiflow. Let $c=c_k$ be the unique generalized critical value of
$V$ in $[a,b]$. Then \be\label{e:5.10} H_q\(V_b,V_a\)\cong H_q
\(V_c,V_c\setminus M_k\),\Hs q=0,1,2,\cdots. \ee
\end{proposition}

\noindent{\bf Proof.} Consider the commutative diagram:
$$
\ba{ccccccccc} H_q(V_a) & \stackrel{i_*}{\longrightarrow}& H_q(V_b)&
\stackrel{j_*}{\longrightarrow}&
 H_q\(V_b,V_a\)& \stackrel{\pa}{\longrightarrow}& H_{q-1}(V_a)& \stackrel{i_*}{\longrightarrow}&
 H_{q-1}(V_b)\\[1ex]
\downarrow{{\tiny i}_* }& &\downarrow{r_* }& &\downarrow{j_*\circ r_* }& &\downarrow{i_* }& &\downarrow{r_*
}\\[1ex]
H_q(V_c\setminus M_k) & \stackrel{i_*}{\longrightarrow}& H_q(V_c)
 & \stackrel{j_*}{\longrightarrow}& H_q\(V_{c},V_{c}\setminus M_k\) &
 \stackrel{\pa}{\longrightarrow}&
H_{q-1}(V_{c}\setminus M_k)  & \stackrel{i_*}{\longrightarrow}&
H_{q-1}(V_{c})
 \ea
$$
The upper and lower rows present the exact homology sequences for
the pairs  $(V_b,V_a)$ and $(V_c,V_c\setminus M_k)$. The
homomorphisms
 $i_*$'s in the vertical arrows are induced by inclusions, and
$r_*$'s  by  deformation retracts from $V_b$ to $V_c$.

As in the proof of the First Deformation Lemma we can show that
$V_a$ is a strong deformation retract of $V_c\setminus M_k$, hence
the vertical arrows number 1 and 4 are isomorphisms. We also infer
from the Second Deformation Lemma that the vertical arrows number 2
and 5 are isomorphisms. Thus the conclusion follows from the
``Five-lemma'' (see \cite{Spa}, Lemma IV.5.11).

\br Assume the hypothesis in Proposition \ref{t:5.1}. Let $U$ be any
neighborhood of $M_k$ with $U\Cap M_i=\ep$ for $i\ne k$. Then by
excision of homologies  we deduce that $$ H_q \(V_c,V_c\setminus
M_k\)\cong H_q \(V_c\cap U,\(V_c\cap U\)\setminus M_k\). $$
Therefore  the concept of critical groups for equilibria of
semiflows coincides with the one for critical points of smooth
functionals \cite{chang2}.

We do not know whether $(\ref{e:5.10}) $ holds true in the  general
case.

\er

The following proposition suggests that the definition of the
critical group of Morse sets is also independent of the choice of
M-L functions.
\begin{proposition}\label{p:5.0}Let $V,L$ be two strict M-L functions of
$\cM$, and let $c=V(M_k)$, $c'=L(M_k)$. Assume  $a<c<b$,
$\alpha<c'<\b$ are such that $c$ and $c'$ are the unique generalized
critical values of $V$ and $L$ in $[a,b]$ and $[\alpha,\b]$,
respectively. Then \be\label{e:5.11} H_q\(V_b,V_a\)\cong H_q
\(L_\b,L_\alpha\),\Hs q=0,1,\cdots.\ee

\end{proposition}

\noindent{\bf Proof.} We may assume that
$$c'=L(M_k)=V(M_k)=c.$$ Otherwise one can replace $L$ by $\Phi=L-c'+c$
(note that $L_\alpha=\Phi_{\alpha-c'+c}$, and
$L_\b=\Phi_{\b-c'+c}$). Define a strict M-L function $F$ of $\cM$
as:
$$F(x)=\max\(V(x),\,W(x)\),\Hs x\in \W.$$
Take an $\ve>0$ sufficient small so that $c$ is the unique
generalized critical value of $F$ in the interval $(c-\ve,c+\ve)$
with
$
(c-\ve,c+\ve)\subset (a,b)\cap(\alpha,\beta).
$
Then \be\label{e:5.13} H_q\(V_b,V_a\)\cong
H_q\(V_{c+\ve},V_{c-\ve}\),\Hs H_q\(L_\b,L_\alpha\)\cong
H_q\(L_{c+\ve},L_{c-\ve}\).\ee

Note that $F_{c\pm\ve}\subset V_{c\pm \ve}\cap L_{c\pm\ve}$. By
sightly modifying the  proof of Theorem \ref{t:4.1} one can show
that $F_{c\pm\ve}$ is a strong deformation retract of both
$V_{c\pm\ve}$ and $L_{c\pm\ve}$. Hence the vertical arrows  number
1, 2, 4 and 5  of the following diagram are isomorphisms:
$$
\ba{ccccccccc} H_q(F_{c-\ve}) & \stackrel{i_*}{\longrightarrow}&
H_q(F_{c+\ve})& \stackrel{j_*}{\longrightarrow}&
 H_q\(F_{c+\ve},F_{c-\ve}\)& \stackrel{\pa}{\longrightarrow}& H_{q-1}(F_{c-\ve})& \stackrel{i_*}{\longrightarrow}& H_{q-1}(F_{c+\ve})\\
\downarrow{i_* }& &\downarrow{i_* }& &\downarrow{i_* }& &\downarrow{i_* }& &\downarrow{i_* }\\
H_q(V_{c-\ve}) & \stackrel{i_*}{\longrightarrow}& H_q(V_{c+\ve})
 & \stackrel{j_*}{\longrightarrow}& H_q\(V_{c+\ve},V_{c-\ve}\) &
 \stackrel{\pa}{\longrightarrow}&
H_{q-1}(V_{c-\ve})  & \stackrel{i_*}{\longrightarrow}&
H_{q-1}(V_{c+\ve})
 \ea
$$
It then follows by Five-lemma that the vertical arrow number 3 is an
isomorphism. That is,
$$
H_q\(F_{c+\ve},F_{c-\ve}\)\cong H_q\(V_{c+\ve},V_{c-\ve}\).$$

Similarly we have
$$
H_q\(F_{c+\ve},F_{c-\ve}\)\cong H_q\(L_{c+\ve},L_{c-\ve}\).$$ Now
the conclusion  follows from (\ref{e:5.13}).
 The proof is
complete.

\Vs To compute  critical groups of Morse sets,  by
 definition one needs to find a strict M-L function of the Morse
 decomposition $\cM$. Here we show that critical groups can be
 successfully computed by using any positively invariant
 neighborhoods of attractors $A_k$ in their attraction basins.
This makes the computation of the critical groups easier and more
flexible.

\bt\label{t:5.4} Let $W_k\subset \W(A_k)$ and $W_{k-1}\subset
\W(A_{k-1})$ be any positively invariant neighborhoods of $A_k$ and
$A_{k-1}$, respectively. Then
\be\label{e:5.14}C_q(M_k)=H_q(W_k,\,W_{k-1}),\Hs q=0,1,\cdots.\ee\et

\noindent{\bf Proof.}  By using primitive Lyapunov functions (see
Theorem \ref{t:2.1}) one can find  positively invariant open
neighborhoods $U_k$ of $A_k$ and  $U_{k-1}$ of $A_{k-1}$ such that
$$U_k\subset W_k,\Hs U_{k-1}\subset W_{k-1}\cap U_k.$$ Further we
infer from Remark \ref{r:2.1} that for $i\in\{k-1,\,k\}$ there
exists a nonnegative function $V_i\in C(X)$ such that
\be\label{e:5.15} V_i(x)\equiv 0\,\,(\mb{on }A_i), \Hs V_i(x)\equiv
1\,\,(\mb{on }X\setminus U_i); \ee \be\label{e:v2} D^+V_i(x)\leq
0,\Hs \A\,x\in U_i\setminus A_i.\ee

Let $V$ be a  strict M-L function of $\cM$. Without loss of
generality we can assume
$$V(x)>0,\Hs \A\,x\in\W.$$ (Otherwise one can use
$\widehat{V}(x)=V(x)-V(M_1)+1$ to replace $V$.\,) Take two positive
numbers $\mu$ and $\lam$ with
$$\mu>V(M_{k-1}),\Hs \lam> V(M_k).$$ Define
$$
L(x)=V(x)+\mu V_{k-1}(x)+\lam V_k(M_k),\Hs x\in \W.
$$
Then $L$ is a strict M-L function of $\cM$. We claim that
\be\label{e:5.16} L_\mu\subset U_{k-1}\subset W_{k-1},\Hs
L_{\lam+\mu}\subset W_k. \ee Indeed, if $x\not\in U_{k-1}$, then by
(\ref{e:5.15}) we have $V_{k-1}(x)=1$. Thus
$$
L(x)>\mu V_{k-1}(x)=\mu.
$$
It follows that $L_\mu\subset U_{k-1}$. Similarly one can show that
$L_{\lam+\mu}\subset W_k$.

Noticing that $$M_{k-1}\subset A_{k-1}\subset A_k,\hs M_k\subset
A_k\-U_{k-1},\hs M_{k+1}\subset X\-U_k\subset X\-U_{k-1},$$ by
(\ref{e:5.15}) we deduce that $$ L(M_{k-1})=V(M_{k-1})+\mu
V_{k-1}(M_{k-1})+\lam V_k(M_{k-1})=V(M_{k-1})<\mu,$$
$$
L(M_{k})=V(M_{k})+\mu V_{k-1}(M_{k})+\lam V_k(M_{k})=V(M_{k})+\mu,$$
$$
L(M_{k+1})=V(M_{k+1})+\mu V_{k-1}(M_{k+1})+\lam
V_k(M_{k+1})=V(M_{k+1})+\mu+\lam>\mu+\lam\,.$$ As $0<V(M_{k})<\lam$,
we see that
$$
\mu<L(M_{k})=V(M_{k})+\mu<\mu+\lam.
$$
Therefore $c=L(M_k)$ is the unique generalized critical value of $L$
in the interval $[\mu,\mu+\lam]$. It follows that
$$C_q(M_k)=H_q(L_{\mu+\lam},\,L_\mu).
$$

On the other hand, since  $W_{k-1}\subset \W(A_{k-1})$ and
$W_k\subset \W(A_k)$ are positively invariant, by slightly modifying
the proof of the First Deformation Lemma one can show that $L_\mu$
and $L_{\mu+\lam}$ are strong deformation retracts of $W_{k-1}$ and
$W_k$, respectively. Using a similar argument as in the proof of
Proposition \ref{p:5.0} below (\ref{e:5.13}), we can easily prove
that
$$
C_q(M_k)=H_q(L_{\mu+\lam},\,L_\mu)\cong H_q(W_k,\,W_{k-1}).
$$

The proof is complete.

\subsection{Morse inequalities}
Now we try to establish Morse inequalities and Morse equations for
attractors. Let \be\label{e:5.3}
m_q=\sum_{k=1}^l\mb{rank}\,C_q\(M_k\),\Hs q=0,1,\cdots. \ee $m_q$ is
called the $q-$th {\bf Morse type number } of $\cM$.

Let $V$ be a given  strict M-L functionof $\cM$,  $c_k=V(M_k)$
($1\leq k\leq l$). Take
 $a,b\in\R$ be such that $$a<c_1<c_2<\cdots<c_l<b.$$ Then $\ep=V_a\subset\a\subset V_b$. Define \be\label{e:5.4}
\beta_q=\beta_q(a,b)=\mb{rank}
\,H_q(V_b,V_a)=\mb{rank}\,H_q(V_b).\ee \bt\label{t:5.2} $(${\bf
Morse inequality}$)$ Suppose that  all the critical groups of each
Morse set $M_k$ are of finite rank. Then for any $q\geq 0$, we have
\be\label{e:5.5} m_q-m_{q-1}+\cdots+(-1)^qm_0\geq
\beta_q-\beta_{q-1}+\cdots+(-1)^q\beta_0.\ee Moreover,
\be\label{e:5.6}
\sum_{q=0}^\8(-1)^qm_q=\sum_{q=0}^\8(-1)^q\beta_q,\ee provided that
the left-hand side of the above equation is convergent.

\et \br\label{r:5.1} Define formal Poincar$\acute{e}$-polynomials
$$
P_\a(t)=\sum_{q=0}^\8 \beta_qt^q,\Hs M_\a(t)=\sum_{q=0}^\8 m_qt^q.
$$
Then $(\ref{e:5.5}) $ can be reformulated in a very simplified
manner: \be\label{e:5.7}M_\a(t)-P_\a(t)=(1+t)Q_\a(t),\ee where
$Q_\a(t)=\sum_{q=0}^\8 \gamma_q\,t^q$ is a a formal polynomial with
$\gamma_q$ being nonnegative integers.

\er

To prove Theorem \ref{t:5.2}, we first need to recall  some basic
facts.

A real function $\Phi$ defined on a suitable family $D(\Phi)$ of
pairs of spaces is said to be {\bf subadditive}, if $Z\subset
Y\subset X$ implies
$$
\Phi(X,Z)\leq \Phi(X,Y)+\Phi(Y,Z).
$$
 If $\Phi$ is subadditive, then for any  $X_0\subset X_1\subset\cdots\subset X_n$ with  $(X_k,X_{k-1})\in D(\Phi)$,
 $$
 \Phi(X_n,X_0)\leq \sum_{k=1}^n\Phi(X_k,X_{k-1}).$$

For any  pair $(X,Y)$ of spaces, set
$$
R_q(X,Y)=\mb{rank}\, H_q(X,Y)\hs(q\mb{-th {\bf Betti number}}).
$$
Define
$$
\Phi_q(X,Y)=\sum_{j=0}^q(-1)^{q-j}R_j(X,Y), \Hs
\chi(X,Y)=\sum_{q=0}^\8(-1)^qR_q(X,Y).
$$
$\chi(X,Y)$ is usually called the {\bf Euler number} of $(X,Y)$.

\bl\label{l:5.1} $\cite{chang,Mil}$  The functions $R_q,\,\Phi_q$
are subadditive, and $\chi$ are additive.\el

\noindent{\bf Proof of Theorem \ref{t:5.2}.} Taking
$$
a=a_0<c_1<a_1<c_2<a_2<\cdots<c_l<a_l=b,$$ by Lemma \ref{l:5.1} one
immediately deduces that
$$
\sum_{i=1}^l\sum_{j=0}^q(-1)^{q-j}R_j\(V_{a_i},\,V_{a_{i-1}}\)\geq\sum_{j=0}^q(-1)^{q-j}R_j\(V_{a_l},\,V_{a_0}\),
$$
that is,
$$
\sum_{j=0}^q(-1)^{q-j}m_j\geq \sum_{j=0}^q(-1)^{q-j}\beta_j.$$

In case $\sum_{q=0}^\8(-1)^qm_q$ is convergent, there is a $q_0$
such that $m_q=0$ for all $q\geq q_0$. It then follows by
(\ref{e:5.5}) that $\beta_q=0$ for $q\geq q_0$, and the conclusion
(\ref{e:5.6}) automatically follows. The proof is complete.

\br As in Theorem \ref{t:5.4}, we can show that for any $b>c_l$ and
 positively invariant neighborhood $W$ of
the attractor $\a$,$$H_*(V_b)\cong H_*(W).$$ Thus
$\beta_q=\mb{rank}\, H_q(V_b)=\mb{rank}\, H_q(W)$. Therefore taking
$W=\W$, one sees that
 $\beta_q$ is the $q$-th Betti number of the attraction basin
$\W$.

\er

 \br\label{r:5.3} If  $\a$ is the global
 attractor of the flow, then we see that $\b_q$ is precisely the $q$-th Betti number
  of the phase space  $X$.
In such a case, it is interesting to note that the right-hand sides
of $(\ref{e:5.5})$ and $(\ref{e:5.6})$ are independent of the
attractor and the flow. It suggests  that  the quantity
$$\mathfrak{M}=\sum_{q=0}^\8(-1)^qm_q=\chi(X)$$ is actually an invariant for dissipative
systems.

A very particular but important case is that $X=M^n$ is an
$n$-dimensional compact $C^1$-manifold, in which all the
 critical groups are of finite rank. We infer from the above argument that for
any Morse decomposition $\cM$ of $M^n$ induced by any flow $S(t)$,
$$
m_0\geq \beta_0,$$
$$
m_1-m_0\geq \beta_1-\beta_0,
$$
$$
\cdots\cdots
$$
$$
m_n-m_{n-1}+\cdots+(-1)^nm_0=\beta_n-\beta_{n-1}+\cdots+(-1)^n\beta_0=\chi(M^n).
$$
\er

\subsection{Morse theory in the quotient phase space}

Let $\~X$ be
 the quotient phase space introduced  in Section 5, and  $\~V([x])=V(x)$ for $[x]\in\~X$.
 We  define the {\bf quotient critical
group} $\~C_*(M_k)$ of $M_k$ to be the homology theory given by
$$
\~C_q(M_k)=H_q\(\~V_b,\,\~V_a\),\Hs q=0,1,\cdots, $$ where $a$ and
$b$ are two real numbers such that $c_k$ is the unique generalized
critical value of $V$ in $[a,b]$, $\~V_R$ denotes the level sets of
$\~V$ in $\~X$.

Clearly all the conclusions concerning the critical group $C_*(M_k)$
hold true for the quotient one. In particular,
 let
$$
\~m_q=\sum_{k=1}^l\mb{rank}\,\~C_q\(M_k\),\Hs
\~\beta_q=\mb{rank}\,H_q(\~V_b), $$ where $b$ is any number with
$b>c_l$. ($\~m_q$ is called the $q-$th {quotient Morse type number
}of $\cM$.) Then  we have
 \bt $($\mb{\bf Quotient Morse inequality}$)$ Assume that
all the quotient critical groups of the Morse sets are of finite
rank. Then for any $q\geq 0$, we have $$
\~m_q-\~m_{q-1}+\cdots+(-1)^q\~m_0\geq
\~\beta_q-\~\beta_{q-1}+\cdots+(-1)^q\~\beta_0.$$ Moreover,
$$
\sum_{q=0}^\8(-1)^q\~m_q=\sum_{q=0}^\8(-1)^q\~\beta_q,$$provided
that the left-hand side of the above equation is convergent. \et

It should be pointed that  in general the critical groups and
quotient critical groups of Morse sets can be different, as is shown
in the following easy example. We suspect that the quotient critical
groups might contain some information of the flow lost by the
critical ones.

\Vs

\noindent {\bf{\em Example }1.} Consider the planar system which
takes the form \be\label{5.e1}\dot{r}=-(r-1)^2,\Hs \dot{\theta}=1\ee
in the polar coordinates. The system has a global attractor
$\a=\ol\cB(0,1)$ with the Morse decomposition $\cM=\{M_1,\,M_2\}$,
where $M_1=0$ and $M_2=S^1$; see Fig 1. \vs

\includegraphics[width=10cm]{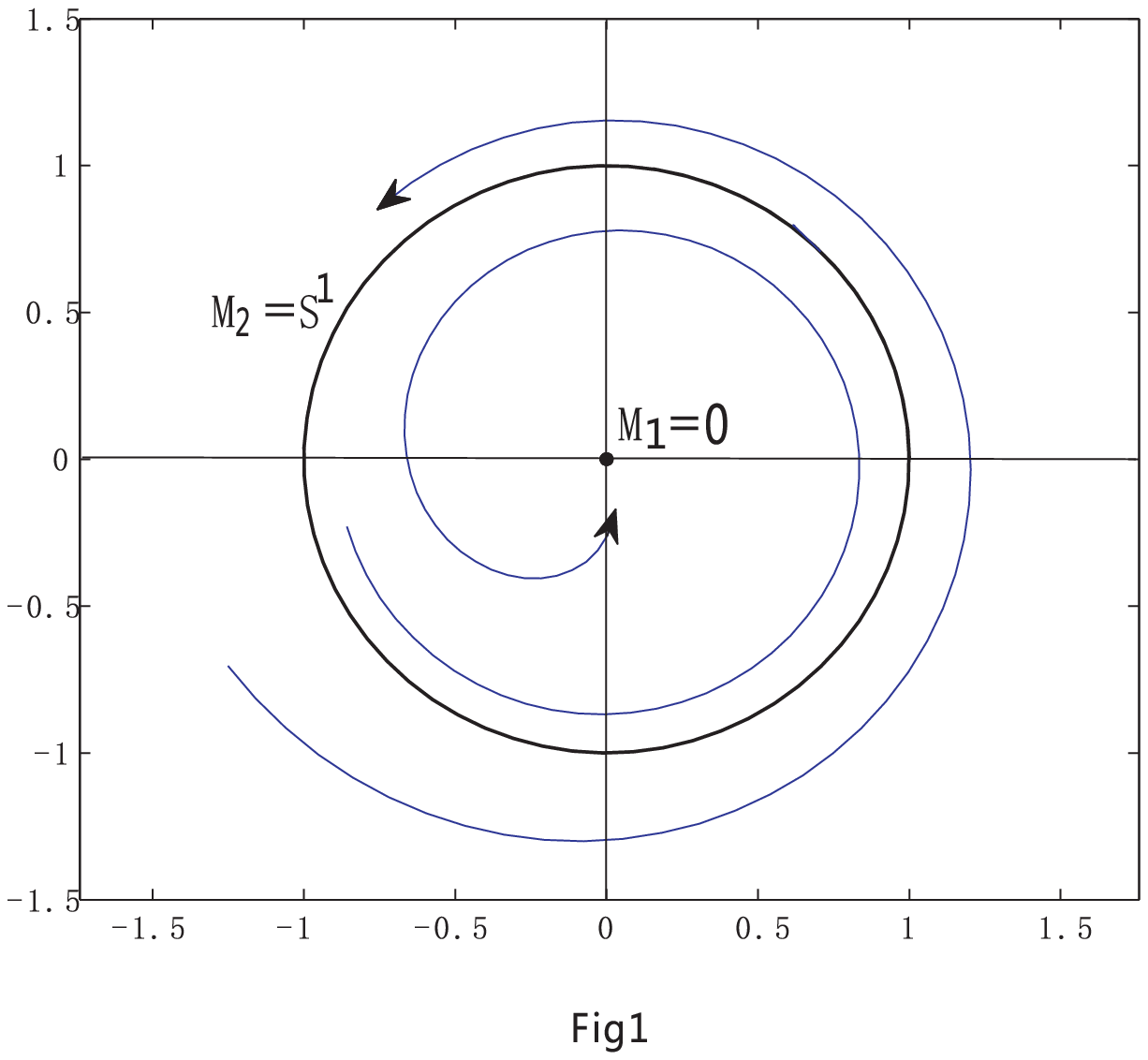}

It is clear that $V(x,y)=r:=\sqrt{x^2+y^2}$ is a Morse-Lyapunov
function of $\cM$. Now let us first compute $\~C_*(M_2)$. Observe
that $ \~V_1=V_1/M_2=S^2$. Since $[0]$ is a strong deformation
retract of $\~V_1\setminus [M_2]$, we deduce by Proposition
\ref{t:5.1} that
$$
\~C_*(M_2)=H_*(\~V_1,\~V_1\setminus[M_2])\cong H_*(S^2,[0]).$$ Using
the  long exact sequence
$$
\cdots\longrightarrow
H_q([0])\stackrel{i_*}{\longrightarrow}H_q(S^2)\stackrel{j_*}{\longrightarrow}H_q(S^2,[0])
\stackrel{\pa}{\longrightarrow} H_{q-1}([0])\longrightarrow\cdots,
$$
one finds that \be\label{e:5.17} H_q(S^2,[0])\cong
H_q(S^2)/\mb{Ker}\,(j_*)=H_q(S^2)/\mb{Im}\,(i_*). \ee Since
$$
H_q(S^2)=\left\{\ba{ll}\mathscr{G},\hs\hs &q=0,2;\\[1ex]
0,&q\ne 0, 2,\ea\right.
$$
we immediately obtain by (\ref{e:5.17}) that
$$
\~C_q(M_2)=\left\{\ba{ll}\mathscr{G},\hs\hs &q=2;\\[1ex]
0,&q\ne 2.\ea\right.
$$

For $C_*(M_2)$ we have by definition that
$$
C_*(M_2)=H_*\(\ol\cB(0,2),\,\ol\cB(0,1/2)\).$$ Since $\ol\cB(0,1/2)$
is a strong deformation retract of $\ol\cB(0,2)$, we see that
$$
C_q(M_2)\cong H_q\(\ol\cB(0,1/2),\,\ol\cB(0,1/2)\)=0,\Hs
q=0,1,2\cdots.$$

\end{document}